\documentclass[11pt,a4paper]{article}
\usepackage{amssymb,latexsym,amsfonts,amsmath,amsthm}

\newtheorem{theorem}{Theorem}[section]
\newtheorem{lemma}[theorem]{Lemma}

\newcommand{\diag}{\mathop{\mathrm{diag}}\nolimits}

\newcommand{\rank}{\mathop{\mathrm{rank}}\nolimits}
\newcommand{\caH}{ \mathop{ \mathcal{H} }\nolimits }
\newcommand{\caS}{ \mathop{ \mathcal{S} }\nolimits }

\newcommand{\caX}{ \mathop{ \mathcal{X} }\nolimits }
\newcommand{\caM}{ \mathop{ \mathcal{M} }\nolimits }
\newcommand{\caN}{ \mathop{ \mathcal{N} }\nolimits }
\newcommand{\caK}{ \mathop{ \mathcal{K} }\nolimits }
\newcommand{\caD}{ \mathop{ \mathcal{D} }\nolimits }

\newcommand{\C}  { \mathop{ \mathbb{C} }\nolimits }
\newcommand{\cl}{\hbox{cl}}

\begin{document}
\title{Pairs of Modules over a Principal Ideal Domain}
\author{Pudji Astuti \thanks{The research of the first author was supported by
Deutscher Akademischer Austauschdienst under Award No.~A/98/25636.}\\
  Departemen Matematika\\
Institut Teknologi Bandung\\
Bandung 40132\\
Indonesia
\and
Harald K. Wimmer\\
Mathematisches Institut\\
Universit\"at W\"urzburg\\
D-97074 W\"urzburg\\
Germany}

\date{}
%\date{\today}

\maketitle

\begin{abstract}
We study pairs of finitely generated modules over a principal ideal domain
and their corresponding matrix representations. We introduce equivalence
relations for
such pairs and determine invariants and canonical forms.

%%%%%%%

\vspace{2cm}
%%%93B25 %systems theory, algebraic methods

\vspace{2cm}

%\vspace{5cm}
%\noindent
%{\bf Address for Correspondence:}\\
%Prof. Dr. H. Wimmer\\
%Mathemat. Institut \\
%Universit\"at W\"urzburg\\
%Am Hubland\\
%D-97074 W\"urzburg\\
%Germany

%\vspace{.2cm}
%\flushleft{
%{\textsf e-mail:}~~\texttt{\small wimmer@mathematik.uni-wuerzburg.de} }\\
%{\textsf Fax:}~~
%+49 931 8\,88\,46\,11\\

\end{abstract}

\section{Introduction}
In this note we study pairs of submodules of $R^n$, where $R$ is a principal
ideal domain. We first review some problems and results on isometries of
$\C ^n$ for which we seek analogies in $R^n$. Let $\caX _1$ and $\caX _2$
be subspaces of $\C ^n$. How are $\caX _1$ and $\caX _2$ situated with
respect to each other? The following theorem, which can be traced back to
Jordan \cite{J},
shows that there exists a suitable orthonormal basis of $\C ^n$ which
displays the respective position of $\caX _1$ and $\caX _2$.
To illustrate  our
purposes it is enough to review the   case where
%no loss of generality to assume
$\dim \caX _1 =
\dim \caX _2 = m$ and $2m = n$. For a subspace $\caX$ of $\C ^n$ let
$P_{\caX}$ denote the orthogonal projection on  $\caX$.

\begin{theorem} \label{Thm.CS}
{\rm{\cite{StS}}}
Assume
\begin{equation} \label{e.X}
\caX _i = X_i \C ^m,\, X_i \in \C ^{n \times m}, \,
X_i^* X_i = I_m, \, i = 1,2, \,\, and \,\, 2m = n.
\end{equation}
Then there are unitary matrices $Q, V_1, V_2,$ such that
\[
        Q X_1 V_1 = \binom{I_m} {0}, \,\,\, Q X_2 V_2 =
 \binom{\Gamma} {\Sigma},
\]
where
\[ \Gamma = \diag (\gamma _1, \dots ,\gamma _m) \quad {\rm{and}} \quad
   \Sigma  = \diag (\sigma _1, \dots ,\sigma _m),
\]
and
\[
     0 \leq \gamma _1 \leq \cdots \leq \gamma _m, \quad
  \sigma _1 \geq \cdots \geq \sigma _m \geq 0,
\]
and
\[
      \gamma _{\mu}^{2} + \sigma _{\mu}^{2} = 1, \,\, \mu = 1,\dots ,m.
\]
The singular values of $P_{\caX _1}(I - P_{\caX _2})$ are $\sigma _1,
\sigma _2,\dots ,\sigma _m,0,\dots ,0$, such that the
numbers $\gamma _{\mu}, \sigma_{\mu}$ are uniquely determined by
$\caX _1$ and $\caX _2$.
\end{theorem}

If $\caX _i, \tilde{\caX _i}, \, i=1,2,$ are subspaces of $\C ^n$   we set
\begin{equation} \label{e.equiv}
(\caX _1, \caX _2) \sim (\tilde{\caX _1},\tilde{\caX _2})
\end{equation}
and we call the two pairs \emph{isometrically equivalent} if there exists
an isometry $\alpha : \C ^{n} \to \C ^{n}$ such that
\begin{equation} \label{e.alpha}
(\tilde{\caX _1},\tilde{\caX _2}) = (\alpha \caX _1, \alpha \caX _2).
\end{equation}
Assume \eqref{e.X} and $\tilde{\caX _i} = \tilde{X _i} \C^{m}, \,
\tilde{X _i}^{*}\tilde{X _i} = I_{m}, \, i = 1,2.$ Then it is obvious that
we have \eqref{e.equiv} if and only if $(X_{1}, X_{2})$ and $(\tilde X_{1},
\tilde X_{2})$ can be transformed into the same canonical form
\[
     \begin{pmatrix} \binom{I}{0}, \, \binom{\Gamma}{\Sigma}
\end{pmatrix}.
\]
That means that one can identify a complete set of invariants under the
equivalence \eqref{e.equiv}.

\begin{theorem} \label{Thm.2}
{\rm{\cite{StS}}} Two pairs  $(\caX _1, \caX _2)$
and  $(\tilde{\caX _1},\tilde{\caX _2})$ are isometrically equivalent if
and only if
\[
       \dim \caX _i = \dim\tilde{\caX _i}, \, i=1,2,
\]
and the singular values of $P_{\caX _1}(I - P_{\caX _2})$ and
${P}_{\tilde{\caX} _1}(I - {P}_{\tilde{\caX}_2})$ are the same.
\end{theorem}

Now let $\caX _i,\tilde{\caX _i}, \,i=1,2,$ be submodules of $R^{n}$.
We say that
the pairs $(\caX _1, \caX _2)$ and  $(\tilde{\caX _1},\tilde{\caX _2})$
are $R-$\emph{unimodular equivalent}, and we write
\[
     (\caX _1, \caX _2) \overset{R}{  \sim  }
(\tilde{\caX _1},\tilde{\caX _2}),
\]
if \eqref{e.alpha} holds for some $R-$automorphism $\alpha : R^{n} \to R^{n}$.
Given two pairs of submodules of $R^{n}$ how can one decide whether they
are $R-$unimodular equivalent? In the case where $\caX _1 + \caX _2$
and $\tilde{\caX _1} + \tilde{\caX _2}$ are direct summands of $R^{n} $
we shall obtain a criterion given in Theorem~\ref{Thm.13} below.
It is known (see e.g. \cite{CR})  that for a
%%%Assuming that $R$ is a
PID one can characterize a \emph{closed}  (or \emph{pure})
submodule $\caX$ of $R^{n}$  by the property that  $\caX$
 is a direct summand of  $R^{n}$ .
If $\caX = X R^{m}$ for some $X \in R^{n \times m}$ then $\caX  $
is closed if and only if all invariant factors of $X$ are $1$.

\begin{theorem}  \label{Thm.13}
 Let $\caX _i,\tilde{\caX _i},\, i=1,2,$ be submodules of $R^{n}$. Assume
that $\caX _1 + \caX _2$ and $\tilde{\caX _1} + \tilde{\caX _2}$
are closed in $R^{n}$. Then the pairs
 $(\caX _1, \caX _2)$ and  $(\tilde{\caX _1},\tilde{\caX _2})$
are $R-$unimodular equivalent if and only if
\begin{equation} \label{e.rank}
  \rank(\caX _1 + \caX _2)  = \rank (\tilde{\caX _1} + \tilde{\caX _2})
\end{equation}
and
\begin{equation} \label{e.quot}
\caX _i \diagup (\caX _{1} \cap \caX _{2}) \cong \tilde{\caX _i}
\diagup \tilde{\caX _1} \cap \tilde{\caX _2}, \, i=1,2,
\end{equation}
hold.
\end{theorem}

Let $(X_1, X_2) \in R^{n \times m_1} \times  R^{n \times m_2}$
and $(\tilde{X}_1, \tilde{X}_2)
 \in R^{n \times \tilde{m}_1} \times  R^{n \times \tilde{m}_2}$
be two pairs of matrices such that $X_i,\tilde{X}_i, \, i=1,2$, have
full column rank. We set
\[
   (X_1, X_2) \overset{u}{\sim} (\tilde{X}_1, \tilde{X}_2)
\]
if
\begin{equation}  \label{e.u}
(\tilde{X}_1, \tilde{X}_2) = (Q X_1 V_1, Q X_2 V_2)
\end{equation}
for some unimodular matrices $Q, V_1, V_2$. If, for $i = 1,2,$ the columns of
$X_i$ and      $\tilde{X}_i$ are bases of $\caX _{i}$ and
$\tilde{\caX _i}$, respectively, then we have
 $ (\caX _1, \caX _2) \overset{R}{  \sim  }
(\tilde{\caX _1},\tilde{\caX _2})$ if and only if
$ (X_1, X_2) \overset{u}{\sim} (\tilde{X}_1, \tilde{X}_2)$. A
`canonical form' of a pair $(X_1, X_2) $ under the transformation
\eqref{e.u} will be needed for the proof of Theorem~\ref{Thm.13}. The
following result is a counterpart of Theorem~\ref{Thm.CS}.

\begin{theorem} \label{Thm.14}
Let $X_1 \in R^{n \times m_1}$   and $X_2 \in R^{n \times m_2}$
have full column rank and let all invariant factors of the matrix
$(X_1 \,\, X_2) \in R^{n \times (m_1 + m_2)}$ be equal to $1$.
Then we have
$(X_1, X_2) \overset{u}{\sim} (Y_1, Y_2)$
where
\begin{equation} \label{e.Y1}
%%%
    Y_1 = \begin{pmatrix} A & 0 \\
                                    0 & I_{m_1 - t}\\
                  0 & 0\\
           0 & 0
    \end{pmatrix}, \, A = \diag (\alpha _1, \dots , \alpha _t), \,\,
 \alpha _t | \cdots | \alpha _1,
%%%
\end{equation}
and
\begin{equation} \label{e.Y2}
Y_2 = \begin{pmatrix} B & 0\\
                               0 & 0\\
                                        0 & I_{m_2 - t}\\
                                      0 & 0
        \end{pmatrix}, \,
                           B = \diag (\beta _1, \dots , \beta _t), \,\,
  \beta _1 | \dots | \beta _t,
\end{equation}
and $0 \leq t \leq \min\{m_1,m_2\}$, and
\begin{equation} \label{e.tfr}
     (\alpha _{\tau}, \beta _{\tau}) = 1, \, \tau = 1, \dots ,t.
\end{equation}
The integer $t$ and the
 elements $\alpha _{\tau}, \beta _{\tau},$ are uniquely determined by
$X_1$ and $X_2$ (up to multiplication by units).
\end{theorem}

The proofs of Theorem~\ref{Thm.13} and Theorem~\ref{Thm.14} will
be given in the next section. We would like to point out that
Theorem 1.4 and the canonical form of pairs $(X_1, X_2)$ under
unimodular equivalence are the main contributions of our paper. For
a general theory of pairs of modules over an arbitrary ring we
refer to \cite{Lv}. In particular, Theorem 1.6 of [3, p.69] on
direct sums of projective rank-1 modules over a commutative ring
goes far beyond our Theorem 1.3.

\section{Proofs}
Let $\caM$ and $\caN$ be  submodules of $R^n$ and $\caM \subseteq
\caN$. The \emph{closure} of $\caM$ in $\caN$  is the submodule
\[\cl(\caM,\caN) = \{ x \in \caN;\, \alpha x \in \caM \,\,
\mathrm{for\,\, some} \,\, \alpha \in R, \, \alpha \neq 0 \}.
\]
%\[
 % \overline{\caM} = \{x \in R^n; \, \alpha x \in \caM \,\,
%\mathrm{for\,\, some}
%\,\, \alpha \in R, \, \alpha \neq 0 \}.
%\]
If $\caN = R^n$ we denote the closure by $\overline{\caM}$.  The
following facts on the closure are known and easy to prove.
%%%(see e.g. \cite{CP}, \cite{CR}).
For two submodules $\caM _1,\caM _2$ of $R^n$ we have
\begin{equation} \label{e.cap}
\overline{\caM _1 \cap \caM _2} = \overline{\caM _1} \cap \overline{\caM _2}
\end{equation}
%%%\end{equation}
and
\begin{equation} \label{e.plus}
\overline{\caM _1 + \caM _2} \supseteq
\overline{\caM _1} + \overline{\caM _2}.
\end{equation}
A submodule $\caM$ is closed in $R^n$ if and only if
$\overline{\caM} = \caM$. From \eqref{e.plus} follows
\begin{equation} \label{e.pluseq}
\caM _1 + \caM _2 =  \overline{\caM _1} + \overline{\caM _2},
\end{equation}
if $\caM _1 + \caM _2 $ is closed.

\begin{lemma} \label{La.21}
Assume that $ \caX _1 + \caX _2$ is a closed submodule of $R^n$.
Define $\caD = \caX _1 \cap \caX _2$. Then
\begin{equation} \label{e.decomp}
\caX _1 + \caX _2 =  \overline{\caD} \oplus \caK _1 \oplus \caK _2
\end{equation}
with
$        \caK _i = \overline{ \caK _i} \subseteq \caX _i, \, i=1,2, $
and
\begin{equation} \label{e.dpl}
\caX _i = (\overline{\caD} \cap \caX _i) \oplus  \caK _i, \, i=1,2.
\end{equation}
%%%
\end{lemma}

\emph{Proof.} Define $\caH_1 = \caX_1 \cap \overline{\caX}_2$ and
$\caH_2 =  {\caX}_2 \cap \overline{\caX}_1 $. Let us show first
that $\caH_i$ is closed in $\caX_i$. Take $i = 1$. If $x$ is in
$\cl(\caH_1,\caX_1)$ then $\alpha x \in \caH_1$ for some nonzero
$\alpha$. Hence $x \in \overline{\caX}_2$ and $x \in \caH_1$. Now
consider $\caS = \caH_1 + \caH_2$. Obviously we have $\caS
\subseteq \overline{\caD}$. To prove the reverse inclusion take
$d \in \overline{\caD},\, d = x_1 + x_2, \ x_i \in \caX_i, $ such
that $\alpha d \in \caD$ for some $\alpha \ne 0.$ Then $\alpha
x_2 \in \caX_1 $ and $x_2 \in \overline{\caX}_2$ and we obtain
$x_2 \in \caH_2$, similarly $x_1 \in \caH_1$, which proves
$\overline{\caD} \subseteq \caS$. Hence we have
\begin{equation}
\label{sumh}
 \caH_1 + \caH_2 = \overline{\caD}\end{equation}
Since, for $i =1, 2,$ the submodule $\caH_i$ is closed in
$\caX_i$ there exists a $\caK_i$ such that $\caX_i = \caH_i
\oplus \caK_i$. From (\ref{sumh}) follows
\begin{equation} \label{e.decomp1}
\caX _1 + \caX _2 =  \overline{\caD} + \caK _1 + \caK _2
\end{equation}
We shall prove the sum in (\ref{e.decomp1}) is direct. Take
 $ y \in  \overline{\caD} \cap (\caK _1 + \caK _2),
 y = k_1 + k_2, \ k_i \in \caK_i$. From $\alpha y \in
\caD, \, \alpha \ne 0$, we obtain $\alpha k_1 \in \caX_1,$ which
yields $k_1 \in \caH_1$. But $\caH_1 \cap \caK_1 = 0$. Hence $k_1
= 0$, and similarly $k_2 = 0$. Therefore $y =0$ and
\begin{equation} \label{e.decomp2}
\caX _1 + \caX _2 =  \overline{\caD} \oplus( \caK _1 + \caK _2).
\end{equation}
It is easy to see that $\caK_1 \cap \caK_2 = 0$, which combined
with (\ref{e.decomp2}) yields (\ref{e.decomp}). To show that
$\caK_i = \overline{\caK}_i$ it suffices to prove that $\caK_i$
is closed in $\caX_1 + \caX_2$. Take $i =1$ and $x \in
\overline{\caK}_1$. If $\alpha x \in \caK_1, \ \alpha x \ne 0$,
then  $x = d + k_1 + k_2, d \in \overline{\caD}, k_i \in \caK_i$,
and (\ref{e.decomp}) imply $\alpha x = \alpha k_1$. Hence $x =
k_1$, and  $\caK_1$ is closed. \hfill $\blacksquare$

%Because of \eqref{e.cap} and \eqref{e.pluseq} we have
%\[
%     \caX _1 + \caX _2 = \overline{ \caX _1 } + \overline{ \caX _2 } =
 %(\overline{ \caX _1 } \cap \overline{ \caX _2 } ) \oplus  \caL =
%  \overline{\caD}  \oplus  \caL
%\]
%for some $\caL$ with $\caL = \overline{\caL }$. Define  $\caK _i =  \caL \cap
%\caX _i, \, i=1,2$. Then $\caL = \caK _1 + \caK _2$. From
%\[
%   (\caL \cap \caX _1) \cap (\caL \cap \caX _2) = \caL \cap \caD
%\subseteq \overline{ \caL \cap \caD} = 0
%\]
%we obtain $\caK _1 \cap \caK _2 = 0$. Hence $ \caX _1 + \caX _2 =
% \overline{\caD} \oplus \caK _1 \oplus \caK _2$ such that
%%
%$R^n =  \overline{\caD} \oplus \caK _1 \oplus \caK _2 \oplus \mathcal{G}$
%%
%for some $\mathcal{G}$. To show that $\caK _i = \overline{ \caK _i}$
%take $i=1$ and $x \in \overline{ \caK _1}$. Then $\alpha x \in \caK _1$
%for some $\alpha \neq 0$. From $x = d + k_1 + k_2 + g, \, d \in
%  \overline{\caD}, k_1 \in \caK _1, k_2 \in \caK _2, g \in \mathcal{G},$
%we deduce $\alpha x = \alpha k_1$, and $x = k_1 \in \caK _1$. Hence
%$\caK _1 = \overline{ \caK _1}$. From \eqref{e.decomp} we immediately
%obtain \eqref{e.dpl}. \hfill $\blacksquare$

{\bf{Proof of Theorem \ref{Thm.14} }}  Let us first consider the case where
$X _1$ and $X _2$ are nonsingular $n \times n$ matrices. The assumption
that all invariant factors of $(X _1 \,\,   X _2)$ are units in $R$ implies
$X_1 S_1 + X_2 S_2 = I_n$ for some $S_1, S_2 \in R^{n \times n}$. Hence
$X _1$ and $X _2$ are left coprime. Put $T = X_1^{-1}X_2$. Then there are
unimodular matrices $V_1$ and  $V_2$ which transform $T$ into
Smith-McMillan form such that
\[
    T = V_1 ^{-1} \diag (\beta _1 / \alpha _1 ,\dots ,\beta _n / \alpha _n) V_2
\]
and
\[
\alpha _n | \cdots | \alpha _1 , \,\,\, \beta _1 | \cdots | \beta _n , \,
\quad \mathrm{and} \quad (\alpha _{\nu}, \beta _{\nu}) = 1,\,
 \nu = 1,\dots , n.
\]
Define
\[
      \hat{X}_1 = \diag(\alpha _1, \dots , \alpha _n) V_1, \,\,
  \hat{X}_2 = \diag (\beta _1, \dots , \beta _n) V_2.
\]
Then $ \hat{X}_1$ and $\hat{X}_2$ are left coprime and
$T = \hat{X}_1 ^{-1 } \hat{X}_2$. Thus we have two left coprime factorizations
of $T$. Hence (see e.g. \cite{V})
$X_1 = Q \hat{X}_1, \,  X_2 = Q \hat{X}_2$ for some unimodular $Q$.

We now deal with the general case and put
$\caX _i = X_i R^{m_i}, \, i=1,2,$  $\caD = \caX _1 \cap \caX _2$, and
$t = \rank \caD$. Let $S$ be unimodular such that the columns of
\begin{equation} \label{e.Smatr}
        S  \left( \begin{matrix}
                             I_t & 0 & 0 \\
        0 & I_{m_1 -t} & 0 \\
                               0 & 0  & I_{m_2 -t} \\
                                                           0 & 0  & 0
       \end{matrix} \right)
\end{equation}
are a basis of  $\mathcal{X}_1 +  \mathcal{X}_2$ which corresponds to the
direct sum in \eqref{e.decomp}. In the following we shall assume
\[
   \mathrm{row \, rank}(X_1 \,\, X_2) = n
\]
or equivalently $\mathcal{X}_1 +  \mathcal{X}_2 = R^n$, which allows us to
discard the bottom row of zero blocks in \eqref{e.Smatr}. Because of
\eqref{e.dpl} a basis of $\mathcal{X}_1$  is given by the columns of a
matrix
       \[
          S \left( \begin{matrix}
       A_1 & 0 \\   0 & I_{m_1 - t} \\ 0  & 0   \end{matrix} \right) \]
%%%%
for some nonsingular $ A_1 \in R^{t \times t}$. Similarly a basis of
 $\mathcal{X}_2$ is of the form
 \[  S \left( \begin{matrix}  B_1  & 0 \\ 0 & 0 \\ 0 &I_{m_2 - t}
\end{matrix} \right) \]
%%%%%
where $B_1$ is nonsingular. Since $\mathcal{X}_1 +  \mathcal{X}_2$   is pure
we have  $A_1 G + B_1 H = I_t$ for some $G, H$. We know already that there
are unimodular matrices $Q_c, W_1, W_2,$ such that $Q_c A_1  W_1 = A$
and $Q_c B_1 W_2 = B$ with $A,B$ as in \eqref{e.Y1} and \eqref{e.Y2}, which
yields the desired canonical form.

To prove uniqueness of $A$ and $B$ in \eqref{e.Y1} and \eqref{e.Y2}  we
note that \eqref{e.tfr} implies that the columns of
$(AB \,\, 0 \,\, 0 \,\, 0)^T$
are a basis of $Q (\mathcal{X}_1 \cap  \mathcal{X}_2)$.
Therefore
\begin{equation} \label{e.q1}
  Q\mathcal{X}_1 \diagup Q (\mathcal{X}_1 \cap  \mathcal{X}_2)
\cong  \mathcal{X}_1 \diagup  (\mathcal{X}_1 \cap  \mathcal{X}_2)
\cong R\diagup \beta _1 R \oplus \cdots \oplus R\diagup \beta _t R,
\end{equation}
%%%
and similarly
\begin{equation} \label{e.q2}
  \mathcal{X}_2 \diagup  (\mathcal{X}_1 \cap  \mathcal{X}_2)
\cong R\diagup \alpha _1 R \oplus \cdots \oplus R\diagup \alpha _t R,
\end{equation}
%%%
which characterizes the entries of $A$ and $B$ in terms of the pairs
$(\mathcal{X}_1 , \mathcal{X}_2) = (X_1 R^{m_1}, X_2 R^{m_2})$.
\hfill $\blacksquare$

{\bf{Proof of Theorem \ref{Thm.13} }} It is obvious that
$    (\caX _1, \caX _2) \overset{R}{  \sim  }
(\tilde{\caX _1},\tilde{\caX _2}) $
implies \eqref{e.rank}   and  \eqref{e.quot}. To prove the converse let
$X_i, \tilde{X}_i, \, i=1,2,$ be basis matrices of $\mathcal{X}_i,
\tilde{\mathcal{X}}_i , \, i=1,2,$ respectively. Then, according to
Theorem~\ref{Thm.14} and \eqref{e.q1} and \eqref{e.q2}, the pairs
$(X_1, X_2)$ and $(\tilde{X}_1, \tilde{X}_2)$ have the same normal form,
which implies
$ (X_1, X_2) \overset{u}{\sim} (\tilde{X}_1, \tilde{X}_2)$ and
$ (\caX _1, \caX _2) \overset{R}{  \sim  }
(\tilde{\caX _1},\tilde{\caX _2}) $. \hfill $\blacksquare$

\vspace*{0.5cm}
 {\bf Acknowledgement.} We are grateful to W.
Schmale  for  valuable comments.

\end{document}